\documentclass[11pt]{article} \usepackage{amssymb}
\usepackage{amsmath}
\usepackage{latexsym, bm}
\usepackage{indentfirst}
\usepackage{CJK}
\usepackage{amsfonts}
\usepackage{amsthm}
\usepackage[all]{xy}
\usepackage[textwidth=5in,textheight=7.5in]{geometry}

\theoremstyle{plain}
\newtheorem{thm}{Theorem}[section]

\newtheorem{lem}[thm]{Lemma}
\newtheorem{cor}[thm]{Corollary}
\newtheorem{prop}[thm]{Proposition}
\theoremstyle{definition}
\newtheorem{defn}[thm]{Definition}

\begin{document}

\title{Odd type generalized complex structures on $4$-manifolds}
\author{Haojie Chen \& Xiaolan Nie}
\date{}
\maketitle
\begin{abstract} We prove that a compact smooth 4-manifold admits generalized complex structures of odd type if and only if it has a transversely
holomorphic $2$-foliation. Consequently, there exist generalized complex structures of odd type on a circle bundle over a closed Seifert fibered $3$-manifold.
\end{abstract}

\section{\textmd{Introduction}}

Generalized complex structure [15, 12]  naturally generalizes symplectic and complex structures, and has many
applications to physical string theory in recent years. In dimension 2, as it reduces to $B$-field transform of symplectic structure or complex structure, the existence problem is trivial. To discuss the existence of generalized complex structures on $4$-manifolds, we divide the problem into two cases: even type and odd type. There are three possibilities in the first case - type zero, type two and type change. While the study of type zero and type two generalized complex structures goes to symplectic topology and complex surface topology, there are surprising results recently on the study of generalized complex structures with type change locus [4,5,10,18,19]. For example, in [4] Calvancanti and Gualtieri proved that there exist type change generalized complex structures on $3\mathbb{C}P^2\#19\overline{\mathbb{C}P^2}$, which does not admit symplectic or complex structures. In [10] and [19] it was proved that the number of components of type change locus can be made arbitrarily large.\let\thefootnote\relax\footnote{Keywords: Generalized complex structures, almost bihermitian structures, transversely holomorphic foliations.} For the odd type case, the type has to be constant type one. In this paper, we study the existence of odd type generalized complex structure on $4$-manifolds. We obtain the following characterization.
\begin{thm} A compact orientable 4-manifold admits odd type generalized
complex structures if and only if it has a structure of
transversely holomorphic 2-foliation.
\end{thm}
Here a transversely holomorphic 2-foliation is a 2-foliation structure whose transition functions are holomorphic (see section 3). Transversely holomorphic 1-foliation on 3-manifold is called
transversely holomorphic flow, with the above theorem, we have the following.
\begin{cor} Let $N^3$ be a compact orientable 3-manifold with a transversely holomorphic
flow. Then there are odd type generalized complex structures
on a circle bundle over $N^3$. \end{cor}
\noindent In particular, the corollary holds for any compact Seifert fibered 3-manifolds, as they have structures of transversely holomorphic flows [3]. This is in analogy with the fact that there are complex structures on $S^1\times N^3$ with $N^3$ a compact Seifert fibered 3-manifold [20].

Using Theorem 1.1, we also prove the following result which first appeared in [2].
\begin{cor} A compact orientable 4-manifold which is diffeomorphic to a
surface bundle over Riemann surface has odd type generalized complex structures.\end{cor}
Actually in [2], Bailey proved the existence of generalized complex structures on a symplectic fiber bundle over a Riemann surface. As surface bundle over surface always admits structures of symplectic fibration, his result implies Corollary 1.3. The method in [2] is to compute certain obstruction classes which differs from ours.\\

The proof of our theorem consists of two parts. First, we prove a general integrability result of generalized complex structures. Given a
generalized complex structure $\mathcal {J}: T\oplus T^*
\longrightarrow T\oplus T^*$ on $M^{2n}$, we show that there always exist two
almost complex structures $j_+, j_-$ (not unique), compatible with
the same Riemannian metric $g$ ( such a triple $(g, j_+, j_-)$ is called an \textbf{almost
bihermitian structure}), and a real two form $b$ such that

$$\mathcal {J} = \frac{1}{2}\left(
                   \begin{array}{cc}
                     1 & 0\\
                     b & 1 \\
                   \end{array}
                 \right)
                 \left(
                   \begin{array}{cc}
                     j_++j_- & \omega_-^{-1}-\omega_+^{-1}\\
                     \omega_+-\omega_- & -(j_+^*+j_-^*) \\
                   \end{array}
                 \right)\left(
                   \begin{array}{cc}
                     1 & 0\\
                     -b & 1 \\
                   \end{array}
                 \right)
$$
where $\omega_\pm(X,Y)=g(j_\pm X, Y)$ for any $X,Y \in T$. Conversely, given
a two form $b$ and an almost bihermitian structure
$(g, j_+, j_-)$, the map $\mathcal {J}$ defined above
will be an almost generalized complex structure. We then derive the integrability conditions of $\mathcal {J}$ in terms of geometric equations of $\omega_+$ and $\omega_-$.

Second, for a compact 4-manifold with transversely holomorphic 2-foliation, we can construct an almost bihermitian structure. To show that it gives an integrable generalized complex structure, we need to verify the above integrability conditions. This is done by choosing the Gauduchon metric in the conformal class of almost Hermitian structures.

The structure of the paper is as follows.

In section 2, we start with preliminaries in generalized complex geometry. After stating the connection between
almost generalized complex structure and almost bihermtian structure, we prove the integrabiligy result. Using it we reprove Gualtieri's theorem on the equivalence between generalized K\"{a}hler structure and the target space geometry of supersymmetric sigma model. When the two almost complex structures from the bihermitian strucuture commute with each other, a much simpler version of the integrability conditions is also achieved.

In section 3, we introduce two definitions of transversely holomorphic foliation and prove their equivalence. After introducing the Gauduchon metric and discussing its existence on almost complex manifolds, we prove Theorem 1.1.\\

\section{\textmd{The almost bihermitian structure and integrability}}
In this paper, all manifolds considered are smooth and orientable. First we give the basic definitions which will be used later. For more details, we refer to [12]. Given a manifold $M$, a natural symmetric pairing on the direct sum of the tangent bundle and cotangent
bundle $T \oplus T^*$ is given by:$$\langle X + \xi , Y + \eta
\rangle = \frac{1}{2}(\eta(X)+ \xi(Y)),
$$ where $X, Y\in T$ and $\xi, \eta\in T^*$. The Courant bracket [13] with respect to a closed real 3-form $H$ of two sections of $T \oplus T^*$ is
defined by: $$[X + \xi, Y + \eta]_H = [X, Y] +\mathcal {L}_X\eta - i_Yd\xi + i_Yi_XH.$$
\theoremstyle{definition}
\begin{defn}An almost generalized complex structure on a manifold
M is an endomorphism $\mathcal{J}$ of $T \oplus T^*$ preserving the natural
pairing such that $\mathcal{J}^2= -id$. A generalized complex structure is an almost generalized complex structure $\mathcal{J}$ such that its $+i$ eigenspace, $L \subset (T
\oplus
T^*)\otimes \mathbb{C}$ is closed under the Courant bracket.\end{defn}
\theoremstyle{plain}
\begin{lem}
Let $\mathcal {J}$  be an almost generalized complex structure on $M^{2n}$, then there exists an almost bihermitian structure  $(g, j_+, j_-)$ and a real two form $b$, such that $$\mathcal {J} = \frac{1}{2}\left(
                   \begin{array}{cc}
                     1 & 0\\
                     b & 1 \\
                   \end{array}
                 \right)
                 \left(
                   \begin{array}{cc}
                     j_++j_- & \omega_-^{-1}-\omega_+^{-1}\\
                     \omega_+-\omega_- & -(j_+^*+j_-^*) \\
                   \end{array}
                 \right)\left(
                   \begin{array}{cc}
                     1 & 0\\
                     -b & 1 \\
                   \end{array}
                 \right)
,$$

where $\omega_\pm=gj_\pm$.
\end{lem}

\emph{Proof}.
As in [12], given an almost generalized complex structure, the structure group of $T\oplus T^*$ can be reduced from $GL(4n)$ into $U(n,n)$ which is homotopic to $U(n) \times U(n)$. Geometrically, we can thus split $T\oplus T^*=C_+\oplus C_-$, where $C_\pm$ are postive/negative definite subbundles with respect to the inner product, and $\mathcal {J}(C_\pm)=C_\pm$. Define $\mathcal {J}_0=-\mathcal {J}$ on $C_+$, and $\mathcal {J}_0=\mathcal {J}$ on $C_-$. Then $(\mathcal {J}, \mathcal {J}_0)$ gives an almost generalized K\"{a}hler structure. By Theorem 2.14 in [14], there always exists
an almost bihermitian structure  $(g, j_+, j_-)$ and a real two form $b$  such that $\mathcal{J}$ is given as above.\qed \vspace{2mm}\\
Now we study the integrability of $\mathcal {J}$.
\begin{thm}
Given an almost generalized
complex structure $\mathcal {J}$ on a smooth manifold $M^{2n}$, let $(g, j_+, j_-,b)$ be an almost bihermitian structure and a real two
form associated to $\mathcal {J}.$ Denote $T_+^{1,0}, T_-^{1,0}\subset T\otimes \mathbb{C}$ to
be the corresponding $+i$-eigenbundles of $j_+, j_-$ and $H$ a closed real three form. Then $\mathcal {J}$ is $H$-integrable if and only if
the following two conditions hold:
\begin{align}
&d\omega_+(X_++X_-, Y_++Y_-,
Z_+)-(\omega_++\omega_-)(X_-,[Y_++Y_-,Z_+])\nonumber\\&+(\omega_++ \omega_-)(Y_-,[X_++X_-,Z_+])-Z_
+\omega_+(X_-,Y_-)\nonumber\\ =&-i(H+db)(X_++X_-,
Y_++Y_-, Z_+)\end{align}
\begin{align}
&d\omega_-(X_++X_-, Y_++Y_-,
Z_-)-(\omega_-+\omega_+)(X_+,[Y_++Y_-,Z_-])\nonumber\\&+(\omega_-+\omega_+)(Y_+,[X_++X_-,Z_-])-Z_-\omega_-(X_+,Y_+)\nonumber\\=&i(H+db)(X_++X_-,
Y_++Y_-, Z_-)
\end{align}
where $X_\pm, Y_\pm, Z_\pm$ are sections of $T_\pm^{1,0}$.
\end{thm}
\emph{Proof}. Given $\mathcal {J}$ as in Lemma $2.2$ and $H$ a closed real 3-form. First consider the $+i$ eigenbundle of
$$\mathcal {J'} = \frac{1}{2}
                 \left(
                   \begin{array}{cc}
                     j_++j_- & \omega_-^{-1}-\omega_+^{-1}\\
                     \omega_+-\omega_- & -(j_+^*+j_-^*) \\
                   \end{array}
                 \right)
.$$
Then direct computation gives $\mathcal {J'}(X + \xi) = i(X + \xi)$ if and only if
$$j_+(X + g^{-1}\xi) = i( X + g^{-1}\xi)\ \ \text{and} \ \ j_-(X - g^{-1}\xi) = i(X - g^{-1}\xi),$$
 i.e. $X + g^{-1}\xi \in T_+^{1,0}$ and $X - g^{-1}\xi \in T_-^{1,0}.$ Note that a vector field $W\in T_+^{1,0}$ if and only if $g(W,Z_+)=0$ for any $Z_+\in T_+^{1,0},$ thus the above conditions give that $(X, \xi)\in +i$ eigenspace of $\mathcal {J'}$ is equivalent to:
\begin{align}
g(X,Z_+)+\xi(Z_+)=g(X,Z_-)-\xi(Z_-)=0
\end{align} for any $Z_\pm\in T_\pm^{1,0}.$ Let $L'$ be the $+i$ eigenspace of $\mathcal {J'}$, and we define the following injective maps\vspace{3mm} \\
$\hspace*{2cm}i_+: T_+^{1,0} \longrightarrow L', \hspace{3cm}i_-:  T_-^{1,0}\longrightarrow L'  \\
 \hspace*{2.6cm}X_+\longmapsto(X_+, gX_+) \hspace{2.6cm}
 X_-\longmapsto (X_-, -gX_-)$.\vspace{3mm}\\
 Then $i_+T_+^{1,0} \cap i_-T_-^{1,0} = 0$ and we have $L' = i_+T_+^{1,0} \oplus i_-T_-^{1,0}$. So for any $(X, \xi) \in L'$, there exists unique $X_+ \in T_+^{1,0}, X_- \in  T_-^{1,0} $ such that $X + \xi = X_+ + X_- + g( X_+ - X_- )$. Now we can get that the $+i$ eigenspace of $\mathcal {J}$
is $$L = e^bL' =  \{X_+ + X_- + b( X_+ + X_- ) + g( X_+ - X_- )| X_\pm\in T_\pm^{1,0}\},$$ then we have
\begin{align*}
&[X_+ + X_- + b(X_+ + X_-) + g(X_+ - X_-), Y_+ + Y_- + b(Y_+ + Y_-)
+ g(Y_+ - Y_-)]_H\\
=& e^b[X_+ + X_- + g(X_+ - X_-), Y_+ + Y_- + g(Y_+ - Y_-)]+i_{Y_++Y_-}i_{X_++X_-}(H+db).
\end{align*}
Since $L'= e^{-b}L$, then\vspace{1.4mm}\\ \hspace*{2cm}$[L,L]_H\subset L\Leftrightarrow [L',L']_{H+db} = e^{-b}[L,L]_H\subset
L'$\vspace{1.4mm}\\
i.e. $L'$ is $H+db$ -involutive.
So we only need to find the integrability condition for $L'$. Let $$X + \xi = X_+ + X_- + g(X_+ - X_-),\
Y + \eta = Y_+ + Y_- + g(Y_+ - Y_-)\in L',$$
the Courant bracket is given by \begin{align*}
[X +\xi, Y + \eta]_{H+db}&= [X_+ + X_-, Y_+ + Y_-] + \mathcal{L}_{X_+ + X_-}g(Y_+ - Y_-)- \\ &i_{Y_+ + Y_-}d(g(X_+ - X_-)) + i_{Y_+ + Y_-}i_{X_+ + X_-}(H+db).
\end{align*}
Use $(3)$ we get that $[X +\xi, Y + \eta]_{H+db}\in L'$
     if and only if
\begin{align}
g&([X_+ + X_-, Y_+ + Y_-], Z_+) + (\mathcal{L}_{X_+ + X_-}g(Y_+ - Y_-)\nonumber \\
&- i_{Y_+ + Y_-}d(g(X_+ - X_-))(Z_+)
=-(H+db)(X_+ + X_-, Y_+ + Y_-, Z_+)
\end{align}
\begin{align}g&([X_+ + X_-, Y_+ + Y_-], Z_-) - (\mathcal{L}_{X_+ +
    X_-}g(Y_+ - Y_-) \nonumber\\
    &- i_{Y_+ + Y_-}d(g(X_+ - X_-))(Z_-)
    = (H+db)(X_+ + X_-, Y_+ + Y_-, Z_-).
    \end{align}
     We use $g = -\omega_\pm j_\pm$ and Cartan formula $\mathcal{L}_X = i_Xd + di_X$ for differential forms to get\\
     \begingroup
\addtolength{\jot}{.4em}
    \begin{align*}&\mathcal{L}_{X_+ + X_-}g(Y_+ - Y_-) - i_{Y_+ + Y_-}d(g(X_+ - X_-))(Z_+) \\
= &i(X_+ + X_-)\omega_-(Y_-, Z_+)
    + i\omega_+(Y_+, [X_+ + X_-, Z_+])
    - i\omega_-(Y_-,[X_+ + X_-, Z_+])\\
&- i(Y_+ + Y_-)\omega_-(X_-, Z_+)
    - iZ_+(\omega_+(X_+, Y_-) - \omega_-(X_-, Y_+))\\
&-i\omega_+(X_+, [Y_+ + Y_-, Z_+]
    + i\omega_-(X_-, [Y_+ + Y_-, Z_+])\\
= &-i(X_+ + X_-)\omega_+(Y_-, Z_+) + i\omega_+(Y_+, [X_+ + X_-, Z_+]) -i\omega_-(Y_-,[X_+ + X_-, Z_+])\\
&+ i(Y_+ + Y_-)\omega_+(X_-, Z_+) - iZ_+(\omega_+(X_+, Y_-) + \omega_+(X_-, Y_+))\\
&-i\omega_+(X_+, [Y_+ + Y_-, Z_+]) + i\omega_-(X_-,
                              [Y_+ + Y_-, Z_+]). \end{align*}
                              \endgroup
 the last equality holds because $\omega_-(Y_-, Z_+) = ig(Y_-, Z_+) = i\omega_+(Y_-, J_+Z_+)
    =\vspace{1.4mm}\\ -\omega_+(Y_-, Z_+). $ Also $g([X_+ + X_-, Y_+ + Y_-], Z_+) = i\omega_+([X_+ + X_-, Y_+ + Y_-], Z_+)$, and $\omega_\pm$ are $(1,1)$ forms with respect to $j_\pm$, then we can rewrite $(4)$ to be
   \begin{align*}
&d\omega_+(X_++X_-, Y_++Y_-,
Z_+)-(\omega_++\omega_-)(X_-,[Y_++Y_-,Z_+])\nonumber\\&+(\omega_++ \omega_-)(Y_-,[X_++X_-,Z_+])-Z_
+\omega_+(X_-,Y_-)\\ =&-i(H+db)(X_++X_-,
Y_++Y_-, Z_+)\end{align*}
     Similarly (5) becomes
   \begin{align*}
&d\omega_-(X_++X_-, Y_++Y_-,
Z_-)-(\omega_-+\omega_+)(X_+,[Y_++Y_-,Z_-])\nonumber\\&+(\omega_-+\omega_+)(Y_+,[X_++X_-,Z_-])-Z_-\omega_-(X_+,Y_+)\\=&i(H+db)(X_++X_-,
Y_++Y_-, Z_-).
\end{align*}
Then the theorem is proved.\qed\vspace{4.5mm}\\
\emph{Remark 1}: The equations $(4)$ and $(5)$ directly indicate that $T_+^{1,0} + T_-^{1,0}$ is involutive. In the case
       $T_+^{1,0} \cap T_-^{1,0} = 0$,  $T_+^{1,0} + T_-^{1,0}=T\otimes\mathbb{C},$ this follows automatically. Otherwise take $Z_+ = Z_- = Z \in T_+^{1,0} \cap T_-^{1,0}$ in $(4), (5)$ and add them to get $g([X_+ + X_-, Y_+ + Y_-], Z) = 0$ for any $Z \in T_+^{1,0} \cap T_-^{1,0}$. Note that the orthogonal\vspace{0.5mm} complement $(T_+^{1,0} \cap T_-^{1,0})^\perp$
is just $T_+^{1,0} + T_-^{1,0} $ as $T_+^{1,0} + T_-^{1,0}\subset(T_+^{1,0} \cap T_-^{1,0})^\perp$ and the two spaces have the same dimensions.  Therefore $[X_+ + X_-, Y_+ + Y_-]\in T_+^{1,0} + T_-^{1,0}$ for any $X_\pm, Y_\pm \in T_\pm^{1,0}$.\vspace{3mm}

Also, we can reprove the following result from Theorem $2.3$, which is originally obtained as Gualtieri's theorem [12] on the equivalence between generalized K\"{a}hler strucuture and target space geometry of $(2,2)$ supersymmetric sigma model.
\begin{prop} (Gualitieri) The almost generalized K$\ddot{a}$hler structure
$$\mathcal {J}_{1,2} = \frac{1}{2}\left(
                   \begin{array}{cc}
                     1 & 0\\
                     b & 1 \\
                   \end{array}
                 \right)
                 \left(
                   \begin{array}{cc}
                     j_+ \pm j_- & -(\omega_+^{-1} \mp \omega_-^{-1})\\
                     \omega_+ \mp \omega_- & -(j_+^* \pm j_-^*) \\
                   \end{array}
                 \right)\left(
                   \begin{array}{cc}
                     1 & 0\\
                     -b & 1 \\
                   \end{array}
                 \right)
,$$ is integrable if and only if
\begin{itemize}
\item $T_\pm^{1,0}$ are involutive.
\item $d_+^c\omega_+ = -d_-^c\omega_- = -(H+db)$
\end{itemize}
\end{prop}
\emph{Proof}. For one direction, consider the almost generalized complex structures $\mathcal{J}_1, \mathcal{J}_2$. The equalities $(4)$ and $(5)$ give integrability of $\mathcal{J}_1$. Since $\mathcal{J}_2$ is obtained from $\mathcal{J}_1$ by replacing $j_-$ by $-j_-$, then $\mathcal{J}_2$ is integrable if and only if the following are satisfied:
\begin{align}
&g([X_+ + \bar{X}_-, Y_+ + \bar{Y}_-], Z_+) + (\mathcal{L}_{X_+ + \bar{X}_-}g(Y_+ - \bar{Y}_-)\nonumber\\
&- i_{Y_+ + \bar{Y}_-}d(g(X_+ - \bar{X}_-))(Z_+)
=-(H+db)(X_+ + \bar{X}_-, Y_+ + \bar{Y}_-, Z_+)\end{align}
\begin{align}&g([X_+ + \bar{X}_-, Y_+ + \bar{Y}_-], \bar{W}_-)- (\mathcal{L}_{X_++
    \bar{X}_-}g(Y_+ - \bar{Y}_-) \nonumber\\
    &- i_{Y_+ + \bar{Y}_-}d(g(X_+ - \bar{X}_-))(\bar{W}_-)
    = (H+db)(X_+ +\bar{X}_-, Y_+ + \bar{Y}_-, \bar{W}_-) \end{align}
where $X_\pm, Y_\pm,Z_\pm \in T_\pm^{1,0}$ and $\bar{X}_-,\bar{Y}_-,\bar{W}_-\in T_-^{0,1}$ denote the complex conjugate.
Let $X_-=Y_-=0$ in $(5)$ and $\bar{X}_-=\bar{Y}_-=0$ in $(7)$, then add them to get:
\begin{align} &g([X_+, Y_+], Z)- (\mathcal{L}_{X_+}g(Y_+)
    - i_{Y_+}d(g(X_+))(Z)\nonumber\\
    &= (H+db)(X_+, Y_+, Z)
\end{align} for any $Z=Z_-+\bar{W}_-$. Since $T_-^{1,0}\oplus T_-^{0,1}=T\otimes \mathbb{C}$, the equation $(8)$ holds for any $Z \in T\otimes \mathbb{C}$. Let $Z=Z_+ \in T_+^{1,0}$ in $(8)$ and $X_-=Y_-=0$ in $(4)$, adding them together gives
$g([X_+, Y_+], Z_+)=0$ for any $Z_+ \in T_+^{1,0}$, i.e. $[X_+, Y_+]\in T_+^{1,0}$ for any $X_+,Y_+\in T_+^{1,0}$. Therefore $T_+^{1,0}$ is involutive. Now $g([X_+, Y_+])= -\omega_+(j_+[X_+, Y_+])=-i\omega_+([X_+, Y_+]),$ we can rewrite the equality $(8)$ to be:\begin{align}i_{[X_+, Y_+]}\omega_+- \mathcal{L}_{X_+}i_{Y_+}\omega_++i_{Y_+}di_{X_+}\omega_+
=i_{Y_+}i_{X_+}i(H+db). \end{align}
Note that $i_{[X_+, Y_+]}=[\mathcal{L}_{X_+},i_{Y_+}]=\mathcal{L}_{X_+}i_{Y_+}-i_{Y_+}\mathcal{L}_{X_+}.$ Use this relation to simplify $(9)$ first and then apply $\mathcal{L}_{X_+}=i_{X_+}d+di_{X_+}$ to get the following:
\begin{align}
i_{Y_+}i_{X_+}d\omega_+=-i_{Y_+}i_{X_+}i(H+db).
\end{align}
Similarly, let $X_+= Y_+=0$ in $(5)$ and take the complex conjugate to get
\begin{align} &g([\bar{X}_-, \bar{Y}_-], \bar{Z}_+) -(\mathcal{L}_{
    \bar{X}_-}g(\bar{Y}_-)
    +i_{\bar{Y}_-}d(g(\bar{X}_-))(\bar{Z}_+)\nonumber\\
    &= -(H+db)(\bar{X}_-, \bar{Y}_-, \bar{Z}_+).
   \end{align} Combine $(6)$ and $(11)$ and use the same argument as above: First we get  for any $Z\in T\otimes \mathbb{C}$.
   \begin{align} &g([\bar{X}_-, \bar{Y}_-], Z) -(\mathcal{L}_{
    \bar{X}_-}g(\bar{Y}_-)
    -i_{\bar{Y}_-}d(g(\bar{X}_-))(Z)\nonumber\\
    &= -(H+db)(\bar{X}_-, \bar{Y}_-, Z).
    \end{align} Then let $Z=\bar{W}_-\in T_-^{0,1}$ in $(12)$ and $X_+=Y_+=0$ in $(7)$ to obtain $g([\bar{X}_-, \bar{Y}_-], \bar{W}_-)=0$ for any $\bar{X}_-, \bar{Y}_-\in T_-^{0,1}$, i.e. $T_-^{0,1}$ is also integrable. Similarly we can use $(12)$ to prove that: \begin{align} i_{\bar{Y}_-}i_{\bar{X}_-}d\omega_-=-i_{\bar{Y}_-}i_{\bar{X}_-}i(H+db). \end{align}
    Note that $H+db$ is a $(1,2)+(2,1)$ form with respect to both $j_\pm$. Use the $(p,q)$ decomposition of forms with respect to $j_\pm$£¬ the equations $(10)$ and $(13)$ are equivalent to:
$$d^c_+\omega+=-d^c_-\omega_-=-(H+db)$$ where $d_{\pm}^c=i(\bar{\partial}_{\pm}-\partial{\pm}).$

For the other direction. If $T_\pm^{1,0}$ are involutive and $d_+^c\omega_+ = -d_-^c\omega_- = -(H+db)$. Consider the integrability of $\mathcal{J}_1$ first. By the above theorem, we only need to show the equalities $(1)$ $(2)$ hold in three cases. If $X_-=Y_-=0$, then $(1)$ is true as $d\omega_+$ and $H+db$ are both $(1,2)+(2,1)$ forms for $j_+$. If $X_+=Y_-=0$, then the left hand side of $(1)$ is equal to \begin{itemize}
\item[]$d\omega_+(X_-,Y_+,Z_+)-(\omega_++\omega_-)(X_-,[Y_+,Z_+])\vspace{1.4mm}\\
=d\omega_+(X_-,Y_+,Z_+)\vspace{1.4mm}\\
= -i(H+db)(X_-, Y_+, Z_+).$
 \end{itemize}
The first equality holds because $T_+^{1,0}$ is involutive and $(\omega_++\omega_-)(X_-,W_+)=0.$ The second equality holds because $d_+^c\omega_+ = -(H+db)$. Now if $X_+=Y_+=0$, as we calculate in theorem $2.2$, the left hand side of $(1)$ becomes
\begin{itemize}

\item[] $d\omega_+(X_-,Y_-,Z_+)-(\omega_++\omega_-)(X_-,[Y_-,Z_+])
+(\omega_-+\omega_+)(Y_-,[X_-,Z_+])\vspace{1.8mm}\\
-Z_+\omega_-(X_-,Y_-)=ig([X_-, Y_-], Z_+) + i((\mathcal{L}_{X_-}g(- Y_-)- i_{Y_-}d(g(- X_-))(Z_+))\vspace{1.8mm}\\
        =-d\omega_-(X_-,Y_-,Z_+)=-i(H+db)(X_-,Y_-,Z_+)$
\end{itemize}
The third equality holds as $d_-^c\omega_- = (H+db)$. So we have proved $(1)$, and similarly for $(2)$ and the integrability of $\mathcal{J}_2$.\qed\vspace{4.5mm}\\
Almost bihermitian structure with two almost complex structures commuting with each other is interesting. With this assumption, we can simplify the equations $(1), (2)$ a lot. To get the simplified version, we need the following observation:

\begin{lem} Let $M$ be a manifold with an almost bihermitian structure $( g, j_+, j_-)$ and $j_+j_- = j_-j_+$. Then the two fundamental forms $\omega_\pm = gj_\pm$ satisfy the following relations:
$$\omega_-(X_+, Y_+) = \omega_+(X_-, Y_-) = 0,\  \omega_+(X, Y) + \omega_-(X, Y) = 0$$ for any $X_\pm, Y_\pm \in T_\pm^{1,0}$ and $X,Y \in T_+^{1,0}+T_-^{1,0}.$
\end{lem}
\emph{Proof}. Since $j_+$ and $j_-$ are both compatible with $g$, also $j_+j_-=j_-j+,$
 we ha\vspace{-3mm}ve \begin{align*}
 \omega_-(X_+, Y_+) &= g(j_-X_+, Y_+) = g(j_+j_-X_+,j_+Y_+)=g(j_-j_+X_+,j_+Y_+)\\&=
 -g(j_-X_+,Y_+)=-\omega_-(X_+, Y_+).\end{align*}
  Thus $\omega_-(X_+, Y_+)=0$.
   Similarly $\omega_+(X_-, Y_-)= 0$. Now we can show the third equality. Given any $X,Y \in T_+^{1,0}+T_-^{1,0}$, there exist $X_\pm, Y_\pm \in T_\pm^{1,0}$ such\vspace{.6mm} that $X=X_+ + X_-, Y = Y_+ + Y_-$.
Note that $\omega_+(X_+, Y_+) = \omega_-(X_-, Y_-) = 0$. Also $\omega_+(X_+, Y_-) = ig(X_+, Y_-) = i\omega_-(X_+, j_-Y_-) = - \omega_-(X_+,Y_-).$ \vspace{.6mm} Therefore $\omega_+(X, Y)+\omega_-(X, Y)=\omega_+(X_-, Y_-)+\omega_-(X_+, Y_+) = 0$.\qed

\begin{cor}
Given an almost generalized
complex structure $\mathcal {J}$ on a smooth manifold $M^{2n}$, let $(g, j_+, j_-, b)$ be an almost bihermitian structure and a real two
form associated to $\mathcal {J}.$ Use the same notation as above. If $j_+$ and $j_-$ commute, then $\mathcal {J}$ is $H$-integrable if and only if $T_+^{1,0}+T_-^{1,0}$ is involutive and the forms $\omega_\pm, b$ and $H$ satisfy the conditions:

$$d\omega_+|_{T_+^{1,0}+T_-^{1,0}}=-d\omega_-|_{T_+^{1,0}+T_-^{1,0}}=-i(H+db)|_{T_+^{1,0}+T_-^{1,0}}$$

\end{cor}
\emph{Proof}. If $\mathcal{J}$ is $H$-integrable, then $T_+^{1,0}+T_-^{1,0}$ is involutive as noted in Remark $1.$ Use the formula $d\omega_\pm(X, Y, Z) = X\omega(Y, Z) + \omega(X, [Y, Z]) + c.p.$ (where $c.p.$\vspace{1.2mm} indicates cyclic permutation) to expand $d\omega_\pm(X_++X_-, Y_++Y_-, Z_+)$. Since $j_+j_-=j_-j_+,$ applying the above lemma, $\omega_+(X, Y) + \omega_-(X, Y) = 0$ for any \vspace{1.3mm} $X,Y \in T_+^{1,0}+T_-^{1,0}$, we obtain that\begin{align} d\omega_+(X_++X_-, Y_++Y_-, Z_+) = - d\omega_-(X_++X_-, Y_++Y_-, Z_+).\end{align}
Similarly we have\begin{align} d\omega_+(X_++X_-, Y_++Y_-, Z_-) = - d\omega_-(X_++X_-, Y_++Y_-, Z_-).\end{align}
Use Lemma $2.5$ again to reduce the interability conditions in Theorem $2.3$ to:
\begin{itemize}
\item[ ] $d\omega_+(X_++X_-, Y_++Y_-, Z_+) = -i(H+db)(X_++X_-, Y_++Y_-, Z_+)$
\item[ ] $d\omega_-(X_++X_-, Y_++Y_-, Z_-) = i(H+db)(X_++X_-, Y_++Y_-, Z_-)$
\end{itemize}
 Hence we get that
\begin{itemize}
\item[ ] $d\omega_+(X_++X_-, Y_++Y_-, Z_++Z_-) = - d\omega_-(X_++X_-, Y_++Y_-, Z_++Z_-)\\= -i(H+db)(X_++X_-, Y_++Y_-, Z_++Z_-).$
\end{itemize}That is, $$d\omega_+|_{T_+^{1,0}+T_-^{1,0}}=-d\omega_-|_{T_+^{1,0}+T_-^{1,0}}=-i(H+db)|_{T_+^{1,0}+T_-^{1,0}}.$$
Conversely, if $T_+^{1,0}+T_-^{1,0}$ is involutive, the above equations and Lemma $2.5$ give $(1), (2)$ in Theorem $2.3$, which imply that $\mathcal{J}$ is $H$-integrable.\qed\vspace{3mm}\\
\emph{Remark 2}: When $j_+=j_-,$ the above conditions indicate $T_+^{1,0}$ is integrable; and when $j_+=-j_-,$ the conditions indicate that $d\omega_+=d\omega_-=0$. In these special cases we see that the generalized complex structures are $b$-transforms of complex or symplectic structures. \vspace{1mm}\\
Another version of the corollary can be obtained from Proposition 3.5 in [2].

\section{\textmd{Proof of Theorem 1.1}}
Let $\mathcal{J}: T\oplus T^*\rightarrow T\oplus T^*$ be an almost generalized complex structure on $M^{2n}$. The type of $\mathcal{J}$ at a point $x$ is defined to be $n-\frac{dim P_xJ(T^*)}{2}$, where $P_x: T_x\oplus T_x^*\rightarrow T_x$ is the projection. It is a fact that the type must be integers and keep the same parity within
a connected component. Therefore, any odd type generalized complex
structure on a 4-manifold must be of constant type one.
The following lemma is from Remark $2.15$ in [14] and Lemma 2 in [1].
\begin{lem}
Given an almost generalized
complex structure $\mathcal {J}$ on a smooth four manifold $M^{4}$, let $(g, j_+, j_-)$ be an associated almost bihermitian structure. Then $\mathcal {J}$ has type one if and only if $j_+, j_-$ induce different orientations. Also, we have $j_+j_-=j_-j_+$ in this case.
\end{lem}
We then get the topological conditions of existence of type one almost generalized complex structures.
\begin{prop} A 4-manifold admits type one almost generalized complex structure if and only if it has an orientable 2-distribution.\end{prop}
\emph{Proof}. By the above lemma, if $(g, j_+, j_-)$ is an associated almost bihermitian structure, then $j_+, j_-$ induce different orientations and $j_+j_-=j_-j_+$. Therefore $ker(j_+\pm j_-)$ will be even dimensional as they are invariant under $j_\pm$. Since $j_+, j_-$ induce different orientations, $ker(j_+\pm j_-)$ each gives an orientable 2-distribution. On the other side, if $E$ is an orientable 2-distribution, choose a Riemannian metric $g$ such that $T=E\oplus N$, where $N$ is orthogonal to $E$. Denote $\rho_E$ to be the rotation of $E$ by $+\pi/2,$  and similarly $\rho_N$ be the rotation of $N$ by $+\pi/2$. Define $j_+=\rho_E+\rho_N, j_-=\rho_E-\rho_N$, then $(g,j_+,j_-)$ is an almost bihermitian structure with $j_+,j_-$ opposite orentations. By Lemma $3.1$ it induces almost generalized complex structure of type one.\qed\vspace{3.5mm}\\
In [16], the topological requirements for the existence of orientable 2-distribution on a compact oriented 4-manifold are given explicitly. For the reader's convenience, we state it here: let $\sigma$ and $\chi$ be the signature and Euler characteristic of the manifold, then it admits an orientable 2-distribution if and only if: (1) if the intersection form is indefinite, then $\sigma + \chi \equiv 0\ mod\ 4$, and $\sigma - \chi \equiv 0\ mod\ 4$; (2) if the intersection form is definite, then $\sigma + \chi \equiv 0\ mod\ 4$, $\sigma - \chi \equiv 0\ mod\ 4$ and $|\sigma| + \chi\geq 0.$\vspace{2mm}\\
Given $\mathcal {J}$ an integrable almost generalized complex structure of type one on $M^4$, by Corollary $2.6$, $T_+^{1,0}+T_-^{1,0}$ is involutive. This turns out to induce structure of transversely holomorphic foliations.
\begin{defn}A $C^\infty$ codimension $p$ foliation $\mathcal{F}$ on
a smooth manifold $M^n$ is given by an open covering $\{U_i\}_{i\in
I}$ of $M$ together with local homeomorphisms $f_i:U_i\rightarrow
\mathbb{R}^n$ such that the
transition functions $f_{ij}$ are smooth and given by:
$f_{ij}(x,y)=(f_{ij}^1(x,y), f_{ij}^2(y))$, where $(x,y)\in \mathbb{R}^{n-p}\times \mathbb{R}^p$, and $f_{ij}^2$ are diffeomorphisms between open sets in $\mathbb{R}^p$.\end{defn}

The subsets $f_i^{-1}(\mathbb{R}^{n-p}\times \{c\})$ in $U_i$ will glue together to give connected immersed submanifolds in $M$, which are called the leaves of $M.$ Use $F$ to denote the tangent distribution of $\mathcal{F}$, whose fiber at each point $x$ is the tangent space to the leaf through $x$. A foliation $\mathcal{F}$ is called \textbf{orientable} if $F$ is an orientable vector bundle. From now on we assume that all the foliations in the rest context are orientable without notification.

\begin{defn}A transversely holomorphic foliation with complex
codimension $m$ is a foliation such that $p=2m$ and $f_{ij}^2$ are
biholomorphisms between open sets in $\mathbb{C}^m$.\end{defn}

Like the complex structure and almost complex structure, we can characterize
the transversely holomorphic foliation in terms of transversely
almost complex structure.\vspace{1.6mm}\\
Let $N=T/F$ be the normal bundle of a foliation $\mathcal{F}$.
\begin{defn}A transversely almost complex
structure of a foliation $\mathcal{F}$ is an endomorphism
$$ j: N\rightarrow N$$ with $j^2=-id$. $(F,j)$ is said to be integrable if $[\pi^{-1}(N^{0,1}),
\pi^{-1}(N^{0,1})]\subset\pi^{-1}(N^{0,1}),$ where $\pi:T\otimes \mathbb{C}\rightarrow N\otimes \mathbb{C}$ is the projection and $N^{0,1}\subset N\otimes \mathbb{C}$ is the $-i$-eigenbundle of $j$.\end{defn}
\begin{prop}{[7]} There is one to one correspondence between transversely holomorphic foliations and foliations with integrable transversely almost complex structure.
\end{prop}
\emph{Proof}: Use the notations as above. For one direction, given a transversely holomorphic foliation $\mathcal{F}$, write the diffemorphisms $f_i$ in Definition $3.2$  to be $(f_i^1, f_i^2),$ where $f_i^2: U_i\longrightarrow \mathbb{C}^m$ are submersions. Then $N|_{U_i}$ can be identified with $f_i^{2*}T\mathbb{C}^m$. If the transition functions $f_{ij}^2$ are holomorphic, then the almost complex structure on $N|_{U_i}$ induced by the holomorphic structure of $T\mathbb{C}^m$ does not depend on $U_i$ and match together to give a transversely almost complex structure. As $\pi^{-1}(N^{0,1})=f_{i*}^{-1}(T\mathbb{R}^{n-2m}\otimes \mathbb{C}\times T^{0,1}\mathbb{C}^m)$ is involutive in each $U_i$, the transversely almost complex structure is integrable.

For the other direction, if $\pi^{-1}(N^{0,1})$ is involutive, since $\pi^{-1}(N^{0,1})+\overline{\pi^{-1}(N^{0,1}})=T\otimes \mathbb{C},$ the existence of transversely holomorphic structure follows from the complex Frobenious theorem [17].\qed\vspace{1.5mm}

We also need to introduce the Gauduchon metric on almost complex manifolds before proving Theorem 1.1.
\begin{defn}Given an almost complex manifold $(M^{2n},j)$, let $T^{*(p,q)}=\wedge^pT^{*(1,0)}\otimes \wedge^qT^{*(0,1)}$ be the bundle with sections $(p,q)$ forms.
Define $\overline{\partial}=\pi_{p,q+1}\circ d: T^{*(p,q)}\rightarrow
T^{*(p,q+1)}$, $\partial=\pi_{p+1,q}\circ d: T^{*(p,q)}\rightarrow
T^{*(p+1,q)}$ to be the formal differential operators. Then an almost Hermitian
metric $g$ is said to be a Gauduchon metric if
$\partial\overline{\partial}\omega^{n-1}=0$ where $\omega=gj$ is the
fundamental form.\end{defn}

The existence of Gauduchon metric is well known on compact complex
manifolds in any conformal class of Hermitian metrics. The
proof actually does not require the integrability of almost complex structures and the existence is also true for almost complex manifolds ([8],[6]).
\begin{prop}{[8]} In any conformal class of almost Hermitian metrics on a compact
connected almost complex manifold, there exists a unique Gauduchon
metric up to scalars.\end{prop}
Now we can prove the theorem.\\
\emph{Proof of Theorem 1.1}: For one direction, if $\mathcal{J}$ is an odd type generalized complex structure on $M^4$, let $(g, j_+, j_-)$ be an associated almost bihermitian structure. From Lemma $3.1$, $j_+, j_-$ induce different orientations and $j_+,j_-$ commute. Therefore $ker(j_+\pm j_-)$ are 2-distributions and we have the decompositions $T_+^{1,0}=A\oplus B$,
$T_-^{1,0}=\bar{A}\oplus B,$ where:\vspace{1.8mm}\\$\hspace*{2.2cm}A=T_+^{1,0}\cap T_-^{0,1}, \hspace{1cm} B=T_+^{1,0}\cap T_-^{1,0}$\vspace{1.8mm}\\are complex line bundles. From Corollary $2.6,$
$T_+^{1,0}+T_-^{1,0}=A\oplus \bar{A}\oplus B$ is involutive. Taking
conjugate gives that $A\oplus \bar{A}\oplus \overline{B}$ is also
involutive. Therefore their intersection $A\oplus \bar{A}=F\otimes
\mathbb{C}$ is involutive,  where $F=ker(j_++j_-)$ is the $2$-distribution.
Moreover, let $N=T/F$, then the projection $\pi: T\otimes \mathbb{C}\rightarrow N\otimes \mathbb{C}$ induces an isomorphism between $B\oplus \bar{B}$ and $N\otimes \mathbb{C}$. Let $N^{0,1}=\pi(B)$, then it defines a transversely almost complex structure
$j: N\rightarrow N$. As \vspace{2mm}\\\hspace*{2.8cm}$\pi^{-1}(N^{0,1})=A\oplus\bar{A}\oplus B=T_+^{1,0}+T_-^{1,0}$\vspace{1.8mm}\\is involutive, $j$ is integrable. Thus it gives a
transversely holomorphic foliation by Proposition $3.6$. Moreover, since $\pi^{-1}(N^{0,1})=T_+^{1,0}+T_-^{1,0}$ is just the projection of the $+i$ eigenbundle of $\mathcal{J}$ to $T\otimes \mathbb{C}$, the transversely holomorphic 2-foliation is uniquely determined by $\mathcal{J}$.

On the other side, given a transversely holomorphic $2$-foliation $(F,
j)$ on a compact 4-manifold $M$, we want to construct an almost
bihermitian structure $(g,j_+,j_-)$ such that $j_+,j_-$ induce different orientations and the conditions in Corollary $2.6$ are satisfied. First embed the bundle $N=T/F$ into $T$ to get a splitting $T=N\oplus F.$ Choose a Riemannian metric $g$ on $M$ such that $g$ is compatible with $j$ on $N$ and $N$ is orthogonal to $F$. Define two almost complex structures on $M$ by $j_+=-j$ on $N$ and $j_+$ to be the rotation by $\pi/2$ on $F$ and\vspace{1.8mm}\\\hspace*{4cm}
$j_-=-j_+|_F+j_+|_N$.\vspace{1.8mm}\\ Then $(g, j_+, j_-)$ is an almost bihermitian structure on $M^4$. By definition, $j_+, j_-$ induce different orientations and $j_+j_-=j_-j_+$. Also we have $T_+^{1,0}+T_-^{1,0}=N^{0,1}\oplus F\otimes \mathbb{C}$. The integrability of $j$ ensures that
$T_+^{1,0}+T_-^{1,0}$ is involutive. By Proposition $3.8$, we can assume
$g$ is a Gauduchon metric with respect to $j_+$ (otherwise, replace it with the Gauduchon metric in its conformal class and it's
still compatible with $j_+, j_-$). Now by Corollary $2.6$, $(g, j_+, j_-)$ gives an $H$-integrable generalized
complex structure if\vspace{1.4mm}\\$\hspace*{1.5cm}d\omega_+|_{T_+^{1,0}+T_-^{1,0}}=-iH|_{T_+^{1,0}+T_-^{1,0}}
=-d\omega_-|_{T_+^{1,0}+T_-^{1,0}}$.\vspace{1.4mm}\\ As $T_+^{1,0}+T_-^{1,0}$ is involutive and $j_+j_-=j_-j_+$, from the equalities $(14)$ and $(15)$ in Corollary $2.6$, we know that
$d\omega_+|_{T_+^{1,0}+T_-^{1,0}}=-d\omega_-|_{T_+^{1,0}+T_-^{1,0}}.$
 Then the theorem follows from the next lemma.\qed
\begin{lem} Given $(g,j_+,j_-)$ an almost bihermitian structure on a four manifold $M^4$ with $j_+, j_-$ inducing opposite orientations
, use the notations as above. If $\partial_+\bar{\partial}_+\omega_+=0$, then there exists a closed real three form $H$ such that $$d\omega_+|_{T_+^{1,0}+T_-^{1,0}}=-iH|_{T_+^{1,0}+T_-^{1,0}}.$$
\end{lem}
\emph{Proof}.  Use the same notions
as above. Let $A=T_+^{1,0}\cap T_-^{0,1}$ and $B=T_+^{1,0}\cap
T_-^{1,0}$ be complex line bundles. As $T\otimes
\mathbb{C}=F\otimes \mathbb{C}\oplus B\oplus \bar{B}$, we have the
splitting $\wedge^3 T^*\otimes \mathbb{C}=\oplus_{k+m+n=3} \wedge^k
\mathbb{F}^*\otimes\wedge^m B^*\otimes \wedge^n\bar{B}^*$, with $\mathbb{F}=F\otimes \mathbb{C}$.
For any three form $\phi\in \Omega^3(M,\mathbb{C})$, we
denote $\phi^{k,m,n}$ to be its components corresponding to the
splitting, where $k\leq2, m,n\leq1$ because of dimension reason. Now to find a close real
three form $H$ such that $d\omega_+|_{T_+^{1,0}+T_-^{1,0}}=-iH|_{T_+^{1,0}+T_-^{1,0}}$, write $H=H^{2,1,0}+H^{2,0,1}+H^{1,1,1}$. Then the equality in the lemma is equivalent to \vspace{1.3mm}\\$H^{2,1,0}=i(d\omega_+)^{2,1,0}$, since $T_+^{1,0}+T_-^{1,0}=F\otimes
\mathbb{C}\oplus B$. As $H$ is real, we get that \vspace{1.4mm}\\ $\hspace*{2.5cm}H^{2,0,1}=\overline{H^{2,1,0}}=-i(d\omega_+)^{2,0,1}$.\\ So only $H^{1,1,1}$ is not determined, and we need to find a three form $H^{1,1,1}$ such that $dH=0$.
Write \vspace{1mm}\\$\hspace*{1.5cm}\bar{\partial}_+\omega_+=
(\bar{\partial}_+\omega_+)^{2,1,0}+(\bar{\partial}_+\omega_+)^{2,0,1}+(\bar{\partial}_+\omega_+)^{1,1,1}.$ \vspace{1.4mm}\\
For any $X\in A\subset T_+^{1,0}, Y\in \bar{A}\subset T_+^{0,1},Z\in B\subset T_+^{1,0}$, we have \vspace{2.5mm}\\ $\hspace*{2cm}(\bar{\partial}_+\omega_+)^{2,1,0}(X,Y,Z)=\bar{\partial}_+\omega_+(X,Y,Z)=0$  \vspace{2.5mm}\\ The second equality holds as $\bar{\partial}_+\omega_+$ is a $(1, 2)$ form with respect to $j_+$. Thus $(\bar{\partial}_+\omega_+)^{2,1,0}$ must be $0$ since it is a section of $\wedge^2(A^*\oplus \bar{A}^*)\otimes B^*.$ Using the fact $d\omega_+=\partial_+\omega_++\bar{\partial}_+\omega_+$ on a 4-manifold, we obtain that \vspace{2.5mm}\\
$\hspace*{1cm} (d\omega_+)^{2,1,0}=(\partial_+\omega_+)^{2,1,0}+(\bar{\partial}_+\omega_+)^{2,1,0}=(\partial_+\omega_+)^{2,1,0},$\vspace{2.5mm}\\ thus $\hspace*{3.5mm}(d\omega_+)^{2,0,1}=
\overline{(d\omega_+)^{2,1,0}}=\overline{(\partial_+\omega_+)^{2,1,0}}=(\bar{\partial}_+\omega_+)^{2,0,1}
$ \hspace{2.5cm} $(16)$\vspace{2.8mm}\\
As $\bar{\partial}_+\omega_+=
(\bar{\partial}_+\omega_+)^{2,0,1}+(\bar{\partial}_+\omega_+)^{1,1,1}$, $\partial_+\bar{\partial}_+\omega_+=0$  gives that \vspace{1.4mm}\\
\hspace*{2cm}
$\partial_+(\bar{\partial}_+\omega_+)^{2,0,1}=-\partial_+(\bar{\partial}_+\omega_+)^{1,1,1} $\vspace{2mm}\\
Similarly $\bar{\partial}_+(\bar{\partial}_+\omega_+)^{2,0,1}=-\bar{\partial}_+(\bar{\partial}_+\omega_+)^{1,1,1}$
as $\bar{\partial}_+\bar{\partial}_+\omega_+=0$. \vspace{2mm}\\ Note that $d\varphi=\partial\varphi+\bar{\partial}\varphi$ for any three form on a $4$-manifold, we actually get \vspace{1.4mm}\\
$\hspace*{3cm}d(\bar{\partial}_+\omega_+)^{1,1,1}= -d(\bar{\partial}_+\omega_+)^{2,0,1}$\\
  Hence by $(16)$, \\$\hspace*{3cm}d(\bar{\partial}_+\omega_+)^{1,1,1}= -d(d\omega_+)^{2,0,1}$\vspace{0.6mm}\\  Now that $dH=0$ if and only if
 \begin{itemize}
 \item[]$\hspace{1.8cm}dH^{1,1,1}=-dH^{2,1,0}-dH^{2,0,1}$
 \item[]$\hspace{3cm}=-id\overline{(d\omega_+)^{2,0,1}}+id(d\omega_+)^{2,0,1}$
 \item[]$\hspace{3cm}=-2$\textmd{Im}$d(d\omega_+)^{2,0,1}$
\item[]$\hspace{3cm}=2$\textmd{Im}$d(\bar{\partial}_+\omega_+)^{1,1,1}$.
 \end{itemize}
 where \textmd{Im} denotes the imaginary part of the form. Let $H^{1,1,1}=2$\textmd{Im}$(\bar{\partial}_+\omega_+)^{1,1,1}$,\vspace{1.3mm}\\
 as $H^{2,0,1}+H^{2,1,0}=-i(d\omega_+)^{2,0,1}+i\overline{(d\omega_+)^{2,0,1}}=2$\textmd{Im}$(d\omega_+)^{2,0,1},$ \vspace{1mm}\\
 \hspace* {2.3cm}$H=2$\textmd{Im}$(\bar{\partial}_+\omega_+)^{1,1,1}+2$\textmd{Im}$(d\omega_+)^{2,0,1} $ \vspace{1.4mm}\\ is a closed real three form such that $d\omega_+|_{T_+^{1,0}+T_-^{1,0}}=-iH|_{T_+^{1,0}+T_-^{1,0}}$.\qed \vspace{2mm}

 The computations in the above lemma are partially helped by those in Bailey's paper. In [2], the obstruction forms are defined (Section 5) for a regular Poisson structure with transverse complex structure. Our lemma can lead to the statement that for a transversely holomorphic 2-foliation on a 4-manifold, the Gauduchon metric $g$ induces a Poisson structure such that the obstruction forms vanish.\vspace{2mm}

\emph{Proof of Corollary 1.2}. By Theorem 1.1, we only need to show there are transversely holomorphic 2-foliations on a circle bundle $X$ over $N^3$ if $N^3$ has a transversely holomorphic flow. Let $\pi: X\rightarrow N^3$ be the projection. Using Definition 3.2, we can find an open covering $\{U_i\}$ of $N^3$ such that the transition function $f_{ij}(x,y)=(f_{ij}^1(x,y), f_{ij}^2(y))$, where $(x,y)\in \mathbb{R}\times \mathbb{C}$ and $f_{ij}^2$ is holomorphic. By replacing $U_i$ by smaller open sets, we can assume that the circle bundle is trivial on $U_i$. Then $\{\pi^{-1}(U_i)\}$ is an open covering of $X$ with the transition functions $g_{ij}(v,x,y)=(g_{ij}^0(v,x,y), f_{ij}^1(x,y), f_{ij}^2(y))$ for $(v,x,y)\in S^1\times \mathbb{R}\times \mathbb{C}$. This defines a transversely holomorphic 2-foliation on $X$. \qed \\

\emph{Proof of Corollary 1.3}. Together with the 2-foliation from the fibering, the complex structure on the base Riemann surface can be pulled back to define a transversely holomorphic 2-foliation structure on the total space. Then the statement follows from Theorem 1.1.\qed \\

Generalized complex structure induces symplectic foliations. The topology of the symplectic leaves is an interesting problem. Using Theorem 1.1 we can construct the following example which gives generalized complex structures of type one on 4-manifold with symplectic leaves noncompact and everywhere dense.\vspace{2mm}\\
\textbf{Example 3.10.} Consider the linear foliation on $T^4=\mathbb{R}^4/\mathbb{Z}^4$. Let $X=(a,b,0,0)$, $Y=(0,0,c,d)$ be two constant vector fields on $\mathbb{R}^4$, with
$a/b, c/d$ irrational numbers, and $E$ be the constant 2-plane field perpendicular to $X, Y$. Then $X, Y$ passes to a 2-foliation $\mathcal{F}$ on $T^4$, and $E$ passes to a two distribution $\tilde{E}$ on $T^4$ which is transversal to $\mathcal{F}$. Any linear transformation $L: E\rightarrow E$ with $L^2=-id$ induces an integrable transversely almost complex structure of $\mathcal{F}$. The leaves of $\mathcal{F}$ are immersed $\mathbb{R}^2$ and are each everywhere dense in $T^4$. By Theorem $1.1$, there are generalized complex structures of type one associated to the linear foliation.\vspace{3mm}

\noindent \emph{Remark 3}: Three-manifolds with orientable transversely holomorphic flows are
classified in [3],[9]. There are six classes of foliations among them,
which are: Seifert fibrations, linear foliations on
$T^3$, strong stable foliations of $T^2$, foliations from holomorphic suspension of $S^2$, affine
foliations on $S^2\times S^1$ and Poincar\'{e} foliations on $S^3$. Recall that the Taubes' conjecture claims that if a closed manifold $S^1\times N^3$ has symplectic structure, then $N^3$ fibers over $S^1$. For a closed 4-manifold $S^1\times N^3$ with transversely holomorphic 2-foliation, it is not known if $N^3$ has transversely holomorphic flows or not.\\

\noindent\textbf{Acknowledgements.}
The authors would like to thank the first author's thesis advisor Tian-Jun Li for many helpful discussions, especially for pointing out Corollary 1.2. They also thank Tedi Draghici for clarifying the Gauduchon metrics and Weiwei Wu for useful conversations. The work in the paper appeared as part of the first author's Ph.D. thesis.

\author{School of Mathematics,\\ University of Minnesota, Minneapolis, MN, 55455}\vspace{2mm}\\
\emph{E-mail address}:  \texttt{chen1512@umn.edu}, \hspace{3mm}\texttt{niexx025@umn.edu}
\end{document}